\documentclass[10pt]{amsart}
\usepackage{amssymb,latexsym,amsthm}

\theoremstyle{plain}
\newtheorem{theorem}{Theorem}
\newtheorem{proposition}[theorem]{Proposition}
\newtheorem{lemma}[theorem]{Lemma}                   
\newtheorem{corollary}[theorem]{Corollary}

\theoremstyle{definition}

\newtheorem{example}[theorem]{Example}

\theoremstyle{remark}
\newtheorem*{remark}{Remark}
\newtheorem*{notation}{Notation}

\DeclareMathOperator{\tr}{tr}

\def\sldoisr{\mathrm{SL}(2,\mathbb{R})}
\def\psldoisr{\mathrm{PSL}(2,\mathbb{R})}
\def\mdoisc{\mathrm{M}(2,\mathbb{C})}

\def\RR{\mathbb{R}}
\def\CC{\mathbb{C}}
\def\ZZ{\mathbb{Z}}
\def\NN{\mathbb{N}}
\def\DD{\mathbb{D}}
\def\SS{\mathbb{S}^1}
\def\eps{\varepsilon}
\def\colon{{:}\;}

\begin{document}

\title[A formula with applications to Lyapunov exponents]
{A formula with some applications to the theory of Lyapunov exponents}
\author[A. Avila and J. Bochi]{Artur Avila and Jairo Bochi}
\date{November 5, 2001}
\thanks{Financial support from Pronex--Dynamical Systems, CNPq 001/2000 and
from Faperj is gratefully acknowledged.}

\address{
Coll\`ege de France -- 3 rue d'Ulm -- 
75005 Paris -- France.}
\email{avila@impa.br}

\address{
IMPA -- Estr. D. Castorina 110 --
22460-320 Rio de Janeiro -- Brazil.}
\email{bochi@impa.br}

\begin{abstract}
We prove an elementary formula about the average expansion of certain
products of $2$ by $2$ matrices.
This permits us to quickly re-obtain an inequality by M. Herman and
a theorem by Dedieu and Shub, both concerning Lyapunov exponents.
Indeed, we show that equality holds in Herman's result.
Finally, we give a result about the growth of the spectral radius of products.
\end{abstract}

\maketitle

\section{Introduction}

A major problem in smooth ergodic theory is to determine whether
a given measure-preserving diffeomorphism has one (or all) Lyapunov
exponents non-zero.
This problem is also of interest in the more general setting of linear cocycles.
However, it is difficult to show the existence of non-zero exponents
without strong conditions like uniform hyperbolicity.

In~\cite{B:Herman}, Herman devised a method to bound the
upper Lyapunov exponent of some cocycles from below and
constructed the first examples of non-uniformly hyperbolic 
two-dimensional systems with a positive exponent.
Such examples are very delicate: it is shown in~\cite{Bochi}
that the exponent of non-hyperbolic cocycles drops to zero with
an arbitrarily small $C^0$-perturbation of the cocycle.

One of the methods of Herman estimates the average of the upper
Lyapunov exponent of systems in a special parametrized family.
While each individual system may be unstable, this average estimate is robust.
Using Herman's estimate, Knill proved in~\cite{B:Knill}
that among bounded measurable $\sldoisr$-cocycles those with a
positive exponent are dense.

This idea -- to consider systems included in some suitable family and to show
that global properties of the family imply
good properties for many individual elements --
is also present in the recent paper~\cite{B:Shub}.
This reasoning has been conjectured to work in more generality
in~\cite {BPSW}.

We will consider the following situation:
take matrices $A_1$,\ldots,$A_n$ in $\sldoisr$.
Of course, the norm $\|A_n \cdots A_1 \|$ can be much smaller
than $\prod\|A_j\|$.
Now we put those matrices inside a family parametrized by a circle:
$A_{j,\theta} = A_j R_\theta$ 
(we indicate by $R_\theta$ a rotation of angle $\theta$).
Instead of looking at the norms, we will deal with the related quantity
$$
N(A)=\log \left( \frac{\|A\|+\|A\|^{-1}}{2} \right ).
$$
In this note we will prove:
$$
\frac{1}{2\pi}\int_0^{2\pi} N(A_{n,\theta} \cdots A_{1,\theta}) \, d\theta
= \sum_{j=1}^{n} N(A_j).
$$
In particular, if the $\|A_j\|$ are large then
$\|A_{n,\theta} \cdots A_{1,\theta} \|$ is of the order of 
$\prod\| A_j \|$ for most values of $\theta$.

The formula allows us to conclude that the mentioned bound of~\cite {B:Herman}
is sharp and also to re-obtain one theorem of~\cite {B:Shub}.

A similar formula, involving spectral radius, also holds.
This motivated us to investigate whether, for cocycles in general,
the spectral radius grows like the norm.
This problem was posed by Cohen in~\cite{B:Cohen}.
The answer, at least in dimension $2$, is no, in general.

\section{The formula}

\begin{notation}
Given a real or complex matrix $A$, we denote:
$$
\|A\|=\sup _{v\neq 0} \frac{\|Av\|}{\|v\|}
\quad\text{where $\|\cdot\|$ is the euclidean norm.}
$$
Also, we denote by $\rho(A)$ the spectral radius, that is, the maximum
absolute value of the eigenvalues of $A$. We have $\rho(A) \leq \|A\|$.
We will indicate by $\sldoisr$ the group of real two-by-two matrices 
with unit determinant and $\psldoisr=\sldoisr / \{\pm I\}$. We define:
$$
N(A)=\log \left( \frac{\|A\|+\|A\|^{-1}}{2} \right )
\quad \text{for $A\in\sldoisr$.}
$$
We define some special matrices in $\sldoisr$:
\begin{align*}
R_{\theta} &=
\begin{pmatrix}
\cos \theta & -\sin \theta \\ 
\sin \theta & \cos \theta
\end{pmatrix}
\quad \text{for $\theta\in \RR$.} \\
H_c &=
\begin{pmatrix}
c & 0      \\
0 & c^{-1}
\end{pmatrix}
\quad \text{for $c\geq 1$.}
\end{align*}
Finally, we indicate by $\DD$ the open unit disk in $\CC$ and by $\SS$ its boundary.
\end{notation}


Our main formula is:
\begin{theorem}\label{form1}
Let $A_1$,\ldots,$A_n \in \sldoisr$. Then
$$
\frac{1}{2\pi}\int_0^{2\pi} 
N( A_n R_{\theta} \cdots A_1 R_{\theta} ) \, d\theta
= \sum_{j=1}^{n} N(A_j).
$$
\end{theorem}

Actually, Theorem~\ref{form1} is a corollary of the formula below:
\begin{theorem}\label{form2}
Let $A_1$,\ldots,$A_n \in \sldoisr$. Then
$$
\frac{1}{2\pi}\int_0^{2\pi} \log\rho
\left ( A_n R_{\theta} \cdots A_1 R_{\theta} \right ) d\theta
= \sum_{j=1}^{n} N(A_j) .
$$
\end{theorem}

Theorems~\ref{form1} and \ref{form2} are proved in
sections~\ref{S:form2} and \ref{S:form1} below.

Notice that $ \log\|A\| - \log 2 < N(A) \leq \log\|A\|$.
Let's give an interpretation of the
quantity $N(A)$ through the following proposition:

\begin{proposition}\label{P:contas}
Let $A \in \sldoisr$. Then
$$
N(A)= \frac{1}{2\pi}\int_0^{2\pi} 
\log \| A(\cos \theta, \sin \theta) \| \, d\theta.$$
\end{proposition}

Therefore the number $N(A)$ can be viewed as
the ``average rate of expansion" of the matrix $A\in\sldoisr$.

\proof
By the polar decomposition theorem, 
one can find numbers $\alpha,\beta \in [0,2\pi]$ and $c \geq 1$
such that $A=R_\beta H_c R_\alpha$.
Moreover, $\|A\| = c$. So we may suppose that $A=H_c$, and we have to prove
$$
\frac{1}{2\pi}\int_0^{2\pi} 
\log \sqrt{c^2\cos^2 \theta+c^{-2}\sin^2\theta} \, d\theta =
\log\left(\frac{c+c^{-1}}{2}\right).
$$

First we calculate
$$
F(b)=\int_0^{\pi} \log (b^2\cos^2 \theta + \sin^2\theta) \, d\theta.
$$
We have
$$
F^{\prime}(b)= 2b \int_0^{\pi} \frac{d\theta}{b^2+\tan^2\theta} 
=2b\int_{-\infty}^{+\infty}\frac{dx}{(b^2+x^2)(1+x^2)}
=\frac{2\pi}{b+1}.
$$
(The last integral can be calculated by residues).
The solution of this differential equation with initial condition
$F(1)=0$ is
$F(b)=2\pi\log\frac{b+1}{2}$.
Therefore
$$
\int_0^{2\pi} 
\log \sqrt{c^2\cos^2 \theta+c^{-2}\sin^2\theta} \, d\theta =
-2\pi\log c + F(c^2) = 2\pi\log\frac{c+c^{-1}}{2}.
$$
\qed

The corollary of Theorem~\ref{form1} below is based on a idea from~\cite{B:Knill}
and justifies the assertion made in the Introduction:
\begin{corollary}
Let $A_1$,\ldots,$A_n \in \sldoisr$ $a>0$ and
$$
E= \left \{ \theta\in [0,2\pi] \colon
\frac{1}{n} \log \|A_n R_\theta \cdots A_1 R_\theta \|
> -a + \frac{1}{n} \sum \log \|A_j\| \right \}.
$$
Let $\nu$ denote the normalized Lebesgue measure in the circle.
Then
$
\nu(E) \geq 1-\frac{\log 2}{a} .
$
\end{corollary}

\proof
Let $f(\theta)=\frac{1}{n} N(A_n R_\theta \cdots A_1 R_\theta)$ and
$M=\frac{1}{n}\sum N(A_j)$.
Let 
$$
F= \{ \theta \colon f(\theta) > M - b \},$$
where $b=a-\log 2$.
It is easy to see that $F \subset E$.
Since $0 \leq f(\theta) \leq M + \log 2$, we have
$$
M =  \int f \, d\nu \leq (M-b)(1-\nu(F)) + (M+\log 2) \nu(F).
$$
This gives
$$
\nu(E)\geq \nu(F) \geq \frac{b}{b+\log 2} = 1- \frac{\log 2}{a}.
$$
\qed

\section{Proof of Theorem~\ref{form2}} \label{S:form2}

The proof is based on complexification methods from~\cite{B:Herman}.

By continuity, we only have to prove the theorem for a dense set of matrices $A_i$.
So we can make the following assumption:
$$
B_\theta = A_n R_\theta \cdots A_1 R_\theta \neq \pm I \quad \text{for all $\theta$.}
$$

Define the following complex matrices:
\begin{align*}
S_z &=
\begin{pmatrix}
\frac{z+z^{-1}}{2}  & -\frac{z-z^{-1}}{2i} \\
\frac{z-z^{-1}}{2i} & \frac{z+z^{-1}}{2}
\end{pmatrix} 
\quad \text{for } z \in \CC^{*}, \\
T_z &=
\begin{pmatrix}
\frac{z^2+1}{2}  & -\frac{z^2-1}{2i} \\
\frac{z^2-1}{2i} & \frac{z^2+1}{2}
\end{pmatrix} 
\quad \text{for } z \in \CC.
\end{align*}
We have $T_z=z S_z$ and $S_{e^{i\theta}}=R_{\theta}$.

Given $A_1$,\ldots,$A_n \in \sldoisr$, we define
$$
C_z = \prod_{j=1}^{n} A_j T_z =
A_n T_z \cdots A_1 T_z \quad \text{for } z\in\CC.
$$

\begin{lemma}\label{L:autoval}
There are holomorphic functions $\lambda_1,\lambda_2 \colon \DD \to \CC$,
which extend continuously to $\overline{\DD}$, such that 
$\{\lambda_1(z),\lambda_2(z)\}$ are the eigenvalues of $C_z$ and
$|\lambda_2(z)|<|\lambda_1(z)|$ for every $z \in \DD$.
\end{lemma}

Lemma~\ref{L:autoval} implies that $\log \rho(C_z)$
is an harmonic function in the disk $\DD$ which extends continuously
to the boundary.
Moreover,
$$
z=e^{i\theta}                                        \quad\Rightarrow\quad
C_z = z^n \prod_{j=1}^{n} A_j R_{\theta}             \quad\Rightarrow\quad
\rho(C_z) = \rho \Big( \prod_{j=1}^{n} A_j R_{\theta} \Big).
$$
Therefore
$$
\frac{1}{2\pi}\int_0^{2\pi} 
\log\rho \Big( \prod_{j=1}^{n} A_j R_{\theta} \Big) \, d\theta 
= \log \rho(C_0).
$$

So the proof of Theorem~\ref {form2} will be complete once we
prove Lemma~\ref{L:autoval} and the lemma below:

\begin{lemma}\label{L:centro}
The eigenvalues of $C_0$ are zero and 
$\prod_{j=1}^{n} \frac{\|A_j\|+\|A_j\|^{-1}}{2}$.
\end{lemma}

\subsection{Proof of Lemma \ref {L:autoval}}

It is enough to show that that the eigenvalues of $C_z$ have different
norms for all $z \in \DD$.
First we obtain the following criteria for identity of their norms:

\begin{lemma}
Let $C \in \mdoisc$ with $\det C \neq 0$ and let
$\lambda_1,\lambda_2$ be the eigenvalues of $C$.
Then $|\lambda_1|=|\lambda_2|$ if and only if
$$
\frac{(\tr C)^2}{4\det C} \in [0,1].
$$
\end{lemma}

\proof
Let $t=\frac{\lambda_1}{\lambda_2}$.
We have $u=\frac{(\tr C)^2}{4\det C} = \frac{1}{4}(t+t^{-1}+2)$.
Then
$   
|\lambda_1|=|\lambda_2| \Leftrightarrow
|t|=1                   \Leftrightarrow
\frac{t+t^{-1}}{2}\in [-1,1]  \Leftrightarrow
u\in[0,1].
$
\qed

We have $\det T_z = z^2$ and so $\det C_z= z^{2n}$.
Therefore $C_z$ has eigenvalues with equal modulus if and only if
$\frac{(\tr C_z)^2}{4z^{2n}} \in [0,1]$, that is, if and only if
$Q(z)=\frac{\tr C_z}{2z^n} \in [-1,1]$.
So to prove lemma~\ref {L:autoval} we must prove that
if $z \in \DD$ then $Q(z)\notin [-1,1]$.
The idea is that since $Q(z)$ is a
rational map of degree at most $2 n$, this can be checked by showing that
the unit circle `exhausts' all preimages of $[-1,1]$.
This we will do with a topological argument.

Let $S = Q^{-1}([-1,1])$.

\begin{lemma} \label {at least}
If $S \cap \SS$ has at least $2 n$ connected components then
$S$ is the union of $2 n$ sub-intervals of $\SS$.
\end{lemma}

\proof
Notice that $\tr C_z$ is a polynomial of degree at most $2 n$, so $Q(z)$ and
$Q'(z)$ are rational maps of degree at most $2 n$.
Since $Q(\SS) \subset \RR$ we know
that there is at least one $0$ of $Q'(z)$ in each component of $\SS
\setminus S$.  In particular there are no zeros of $Q'(z)$ in $S$, which
implies that each connected component of $S \cap \SS$ is mapped
diffeomorphically onto $[-1,1]$.
\qed

Define the following sets:
\begin{align*}
X &= \psldoisr \setminus \{ I \},        \\
Y &= \{\, A \in X \colon |\tr A|\leq 2\,\}.
&& \text{($|\tr|$ is well-defined in $\psldoisr$)}
\end{align*}

\begin{lemma}
There exists a continuous function $F \colon X \to \SS$ such that
$F^{-1}(\{1\})=Y$ and that the induced homeomorphism
$F_\# \colon \pi_1(X) \to \pi_1(\SS)$ is an isomorphism.
\end{lemma}

\proof
Let $A\in\psldoisr$. If $A\in Y$ then
we define $F(A)=1$. Otherwise $A$ has two eigendirections
$\pm v$ and $\pm w$, where $v,w\in\SS\subset \CC$, with associated eigenvalues
$\lambda$ and $\lambda^{-1}$, where $|\lambda|>1$.
We then define $F(A)$ as $v^2/w^2$.
It is easy to see that $F$ is continuous at every $A\neq I$.

We have $\pi_1(X)=\ZZ$ (it is equal to $\pi_1(\psldoisr)$, since $\psldoisr$
is a three-dimensional manifold) and so it is enough to exhibit
a closed path $\gamma$ generating $\pi_1(X)$
such that $F \circ \gamma$ has degree one.

Let $\gamma \colon \SS \to X$ be defined as
$e^{i \theta} \mapsto R_{\theta/2} M$,
where $M \in X$ is symmetric.
Using the identities  $F(A^T)=\overline{F(A)}$ and
$F(R_\theta^{-1}AR_\theta)=F(A)$, we easily
see that $F \circ \gamma$ commutes with conjugacy
($F \circ \gamma(\overline{z})=\overline{F \circ \gamma(z)}$).
Furthermore,
$e^{i \theta}=1$ is the only value such that
$\gamma(e^{i \theta})$ is symmetric, which is equivalent to
$F(\gamma(e^{i \theta}))=-1$.
This implies that $F \circ \gamma$ has degree one.
\qed

Let $B_{\theta}=\prod_{j=1}^{n} A_j R_{\theta}$.
We have $Q(e^{i \theta})=\frac{1}{2} \tr B_{\theta}$.
Notice that $e^{i \theta} \in S$ if and only if $B_\theta$ is an 
elliptic or parabolic matrix in $\sldoisr$.

\begin{lemma} \label{has}
$S\cap\SS$ has at least $2n$ connected components.
\end{lemma}

\proof
Let $g \colon \SS \to \psldoisr$ be defined by $e^{i \theta} \mapsto R_{\theta/2}$. 
This is clearly a generator of the fundamental group of $\psldoisr$.

Notice that $e^{i \theta} \mapsto B_\theta$ can also be seen as a path in
$\psldoisr$ and it follows from the definition that it is homotopic to
$g^{2 n}$ (by the fact that $\psldoisr$ is a group).
Now we use the assumption made at the beginning:
for all $\theta$, $B_\theta \neq \pm I$.

In this case, the path $F \circ B_\theta$ has degree $2 n$ and therefore
the preimage of $1$ has at least $2 n$ connected components.
This set coincides with $S \cap \SS$.
\qed

Lemmas~\ref{at least} and~\ref{has} imply that $S\subset \SS$,
and Lemma~\ref{L:autoval} is proved.
\qed

\subsection{Proof of Lemma \ref {L:centro}}
  
We will list some facts to be used:
\begin{itemize}
\item[(1)]
One can find numbers $\alpha_j,\beta_j \in [0,2\pi],\, c_j \geq 1$
such that $A_j=R_{\beta_j}H_{c_j}R_{\alpha_j}$ for each $j$.
Moreover, $\|A_j\|=c_j$.

\item[(2)]
$A,B\in\sldoisr\quad\Rightarrow\quad\rho(AB)=\rho(BA)$.

\item[(3)]
For every $\theta\in\RR$, $R_{\theta}T_0 = T_0 R_{\theta} =
e^{-i\theta}T_0$.

\end{itemize}

Part (1) is the polar decomposition theorem.
For (2), notice that the spectral radius depends only on the trace.
For (3), we use that $S_z S_w=S_{zw}$.
This implies $T_z T_w=T_{z w}$ and
$$
R_{\theta}T_0 = S_{e^{i\theta}}T_0 = e^{-i\theta}T_{e^{i\theta}} T_0 =
e^{-i\theta} T_0.
$$

Using (1), (2) and (3), we obtain
$$
\rho(C_0) = \rho \Big( \prod_{j=1}^{n}A_j T_0      \Big)
          = \rho \Big( \prod_{j=1}^{n} T_0 H_{c_j} \Big).
$$
Each matrix $T_0 H_{c_j}$ has an eigenvector
$(-i,1)$ with corresponding eigenvalue
$\frac{c_j+c_j^{-1}}{2}$.
Therefore $\prod_{j=1}^{n} T_0 H_{c_j}$ has an eigenvalue
$$
\prod_j (c_j+c_j^{-1})/2,
$$
while $C_0$ is not invertible.
This proves Lemma~\ref {L:centro} and hence Theorem~\ref{form2}.
\qed

\section{Proof of Theorem \ref{form1}}  \label{S:form1}

Let $B_\theta=A_n R_\theta \cdots A_1 R_\theta$.  Then, fixing $\theta$ we have,
by Theorem~\ref{form2},
$$
\frac{1}{2\pi}\int_0^{2\pi} \log\rho
(B_\theta R_{\theta'}) d\theta'=
N(B_\theta).
$$
On the other hand, fixing $\theta'$ we have, again by Theorem~\ref{form2},
\begin{align*}
\frac{1}{2\pi}\int_0^{2\pi} \log\rho
(B_\theta R_{\theta'}) d\theta &=
\frac{1}{2\pi}\int_0^{2\pi} \log\rho
\left ( A_n R_{\theta} \cdots (A_1 R_{\theta'}) R_{\theta} \right )
d\theta\\
&=  N(A_1 R_{\theta'})+ \sum_{j=2}^{n} N(A_j) 
= \sum_{j=1}^{n} N(A_j).
\end{align*}
Then
$$
\frac {1} {2\pi} \int_0^{2\pi}
N(B_\theta) \, d\theta
= \frac {1} {2\pi} \int_0^{2\pi} \frac{1}{2\pi}\int_0^{2\pi} \log\rho
(B_\theta R_{\theta'}) \, d\theta d\theta' \\
= \sum_{j=1}^{n} N(A_j).
$$
This proves Theorem~\ref{form1}.
\qed

\begin{remark}
Inversely, Theorem~\ref{form2} could be quickly deduced from Theorem~\ref{form1}, using
$$\log\rho(A)=\lim_{n\to \infty}\frac{1}{n}\log\|A^n\|=\lim_{n\to \infty}\frac{N(A^n)}{n}.$$
\end{remark}

\section{Herman's inequality re-obtained}
\label{S:Herman}

Let $(X,\mu)$ be a probability space and $T \colon X\to X$ an
ergodic transformation. Let $A \colon X \to \sldoisr$ be a measurable function
satisfying the integrability condition
$$
\int \log \|A\| d\mu < \infty.
$$
We denote for $x\in X$ and $n\in \NN$ ,
$$
A^{n}(x)= A(T^{n-1}x)\cdots A(x) \, .
$$
The function $A$ is called a \textbf{linear cocycle}.
In these conditions,
there exists (see~\cite{B:Fu} or~\cite{B:Led}) a number $\lambda ^{+}(A)\geq 0$,
called the \textbf{upper Lyapunov exponent}, such that 
$$
\lambda ^{+}(A)=\lim_{n\to +\infty }\frac{1}{n}\log \left\|A^{n}(x)\right\| 
\quad \text{for $\mu$-a.e. $x\in X$.}
$$

For $\theta\in\RR$, we define a cocycle $AR_\theta$ by 
$(AR_\theta)(x)=A(x)R_\theta$.
Clearly, $\theta \mapsto \lambda^{+}(AR_\theta)$ is a measurable function.

We now state Herman's inequality:

\begin{theorem}[\cite{B:Herman}, \S 6.2, see also \cite{B:Knill}]
\label{desigualdade}
If $T$, $\mu$ and $A$ are as above then
$$
\frac{1}{2\pi}\int_{0}^{2\pi} \lambda ^{+}(AR_{\theta}) d\theta
\geq \int_{X}\log
\left ( \frac{\|A(x)\|+\|A(x)\|^{-1}}{2} \right ) d\mu(x)\, .
$$
\end{theorem}

\begin{remark}
Herman's inequality was stated in a different (but equivalent) way, involving
the Iwasawa decomposition.
\end{remark}

We will re-obtain Theorem~\ref{desigualdade} and also show that equality holds.
 
\begin{theorem} \label{igualdade}
If $T$, $\mu$ and $A$ are as above then
$$
\frac{1}{2\pi}\int_{0}^{2\pi} \lambda ^{+}(AR_{\theta}) d\theta=
\int_{X}\log
\left ( \frac{\|A(x)\|+\|A(x)\|^{-1}}{2} \right ) d\mu(x)\, .
$$
\end{theorem}

\proof
Recall that $N(A) \leq \log \|A\| < \log 2 + N(A)$.
By Theorem~\ref{form1},
$$
\sum_{j=0}^{n-1} N(A (T^j(x))
\leq
\frac{1}{2\pi}\int_{0}^{2\pi} \log \| (AR_{\theta})^{n}(x) \| \, d\theta
\leq
\log 2 + \sum_{j=0}^{n-1} N(A (T^j(x)).
$$
Therefore, by Birkhoff's theorem,
$$\lim_{n\to\infty}
\frac{1}{2\pi}\int_{0}^{2\pi} \frac{1}{n} \log \| (AR_{\theta})^{n}(x) \| \, d\theta
= \int N(A(x)) \, d\mu(x)
\quad \text{for a.e. $x$.}
$$
To finish the proof we must check that Dominated Convergence applies.
We have
$$
0 \leq
\frac{1}{n} \log \| (AR_\theta)^n (x) \| \leq
\frac{1}{n} \sum_{j=0}^{n-1} \log \| A(T^j x) \|
= f_n (x).
$$
$\{ f_n \}$ is the sequence of Birkhoff means of the function
$\log \|A\| \in L^1 (\mu)$.
In particular, $\{ f_n (x)\}$ is bounded for a.e. $x$.
\qed

\begin{example}
Consider the cocycle (\cite{B:Herman},~\S~4.1) where
$T \colon \SS \to \SS$ is an (uniquely ergodic) irrational rotation, 
$A \colon \SS \to \sldoisr$ is given by 
$A(e^{it})=H_c R_t$ and $c\geq1$ is fixed.
We have 
$(AR_\theta)^n (z) = A^n(e^{i\theta} z)$ and therefore
$\lambda^{+}(A)=\lambda^{+}(AR_\theta)$ for all $\theta$.
It follows from Theorem~\ref{igualdade} that
$\lambda^{+}(A)=\log \left( \frac{c+c^{-1}}{2}\right)$.
\end{example}

\section{A theorem by Dedieu and Shub re-obtained}
\label{S:Shub}

We will use Proposition~\ref{P:contas} and Theorem~\ref{form2} (in the case $n=1$) 
to give another proof of the following theorem by Dedieu and Shub:

\begin{theorem}[\cite{B:Shub}]
Let $\mu$ be a probability measure in $\sldoisr$ 
such that the integral
$\int\log\|A\|\, d\mu(A)$ is finite.
Suppose that $\mu$ is invariant by rotations,
that is, $R_\theta^* \mu=\mu$ for all $\theta$.
Let $A_1$,$A_2$,\ldots$\in\sldoisr$
be independent random matrices with law $\mu$
and consider the associated upper Lyapunov exponent:
$$
\lambda^{+}=\lim_{n\to +\infty }\frac{1}{n}\log \|A_n\cdots A_1 \| 
\quad \text{(w.p. $1$).}
$$
Then
$$
\lambda^{+} = \int_{\sldoisr} \log\rho(A) \, d\mu(A).
$$
\end{theorem}

\proof
\begin{align*}
\lambda^{+} &=
\int \int_{0}^{2\pi} \log \|A e^{i\theta}\| \, \frac{d\theta}{2\pi} \, d\mu(A)
&&\text{(by Furstenberg's formula, see \cite{B:Led})} \\
&= \int N(A) d\mu(A)
&&\text{(by Proposition \ref{P:contas})} \\
&= \int \int_{0}^{2\pi} 
\log\rho(AR_{\theta}) \, \frac{d\theta}{2\pi} \, d\mu(A)
&&\text{(by Theorem \ref{form2} with $n=1$)} \\
&= \int_{0}^{2\pi} \int
\log\rho(AR_{\theta}) \, d\mu(A) \, \frac{d\theta}{2\pi}
&&\text{(since $\log \rho (AR_\theta) \leq \log\|A\| \in L^1$)} \\
&= \int \log\rho(A) \, d\mu(A)
&&\text{(since $\mu$ is invariant by rotations).}
\end{align*}
\qed

\section{Growth of the spectral radius}

Let $X$, $\mu$, $T$ and $A$ be as in section~\ref{S:Herman}.
In view of our results,
it is somewhat natural to ask about the behavior the spectral radius
of the matrix $A^n(x)$ when $n\to\infty$.
This question was already raised in~\cite{B:Cohen}.
We have the following result:

\begin{theorem}  \label{T:expoente espectral}
Suppose $T$ is invertible. Then
for $\mu$-a.e. $x\in X$,
$$
\limsup_{n\to\infty} \frac{1}{n} \log\rho (A^n(x)) = \lambda^{+}(A).
$$
\end{theorem}

Before giving the proof, we point out that in general
the limit of $\frac{1}{n} \log\rho(A^n(x))$
does not exist.
Furthermore, the relation 
\begin{equation}
\limsup_{n\to\infty} \frac{1}{n} \int \log\rho(A^n) \, d\mu = \lambda^{+}(A)
\tag{$\ast $}
\end{equation}
is in general \emph{false}, as is shown by the following:
\begin{example}
Let $X=\{0,1\}^{\ZZ}$, $\mu$ be the $(\frac12, \frac12)$-Bernoulli measure and
let $T \colon X \to X$ be the left shift.
We define a cocycle $A \colon X \to \sldoisr$ by:
$$
A(\{x_i\}_{i\in \ZZ})=
\begin{cases}
H = H_2         &\text{if $x_0=1$,} \\
I               &\text{if $(x_{-1},x_0,x_{1})=(0,0,0)$ or $(1,0,1)$,} \\
R = R_{\pi/2}   &\text{if $(x_{-1},x_0,x_{1})=(1,0,0)$ or $(0,0,1)$.}
\end{cases}
$$
Given any  sequence $x = \{x_i\}_{i\in \ZZ}$, split it in minimal blocks
starting with $1$, as for instance,
$$
\dots (10) (1) (1000) (100) (10) (100000) (1) \dots
$$
The corresponding splitting for the sequence $\{A(T^i(x))\}$ is,
in this case,
$$
\dots (HI) (H) (HRIR) (HRR) (HI) (HRIIIR) (H) \dots
$$
The product of the matrices in each block is always $\pm H$.
It follows that $ \lambda^{+}(A) = \frac{\log 2}{2}$.
On the other hand, making substitutions $R^2 = -I$
in the product $A^n(x)$, we obtain one of the possibilities:
$\pm H^k$, $\pm H^k R$, $\pm R H^k$ or $\pm R H^k R$.
Since $\rho(H^k R) = \rho(R H^k) = 1$, we have
$\rho (A^n(x)) = 1$ infinitely often for a.e. $x$.
Besides, it's not hard to show that~($\ast$) does not hold.
\end{example}

\proof[Proof of Theorem~\ref{T:expoente espectral}]
We may regard the problem as being posed in $\psldoisr$ instead of $\sldoisr$.
Suppose that $\lambda^{+}(A)>0$ (otherwise there is nothing to prove).
Consider (see~\cite{B:Led}) the Oseledets splitting
$\RR^2= E^{+}(x) \oplus E^{-}(x)$,
defined for a.e. $x \in X$, where $E^{+}$ (resp. $E^{-}$) is associated to
the exponent $\lambda^{+}(A)$ (resp. $-\lambda^{+}(A)$).
By Oseledets' theorem,
$$
\lim_{n\to\infty} \frac{1}{n}
\log\sin\measuredangle\left(E^{+}(T^n x),E^{-}(T^n x)\right) = 0.
$$

For each $x$, take $B(x)\in \psldoisr$ that sends the direction $\RR (1,0)$
(resp. $\RR (0,1)$) to the direction $E^{+}(x)$ (resp. $E^{-}(x)$).
This defines a{.}e{.} a measurable function $B \colon X \to \psldoisr$
such that
$$
\lim_{n\to\infty} \frac{1}{n} \log \| B(T^n x) \| = 0
\quad \text{for a.e. $x$.}
$$

We claim that
$$
\liminf_{n\to\infty} \|B(x)^{-1} B(T^n x) -I \| = 0
\quad \text{for a.e. $x$.}
$$
To prove it, let $\eps>0$.
Consider a countable cover of $\psldoisr$ by open sets
$$
U_j= \{ \, M\in \psldoisr \colon
\| M-M_j\| < \delta_j \, \},
\quad \text{where $2 \delta_j (\|M_j\| + \delta_j) < \eps$.}
$$
Define $V_j = B^{-1}(U_j)\subset X$ and
$$
\tilde{V_j} = \{ \, x\in V_j \colon
T^n(x) \in V_j \text{ for infinitely many } n \in \NN \, \}.
$$
By Poincar\'{e}'s recurrence theorem, $\mu(\tilde{V_j})=\mu(V_j)$.
If $x\in \tilde{V_j}$ then, for infinitely many $n \in \NN$, we have
$$
\| B(x)^{-1} B(T^n x) - I \| \leq \| B(T^n x) - B(x) \| \cdot \|B(x)\|
< 2\delta_j (\|M_j\| + \delta_j) < \eps.
$$
Therefore $\liminf \|B(x)^{-1} B(T^n x) - I \| \leq \eps$ for every $x$ in the
full measure set $\bigcup \tilde{V_j}$.
This proves the claim.

\smallskip

To prove the Theorem it's enough 
(since $\rho + \rho^{-1} = \max \, \{ | \tr |, \, 2 \}$)
to show that
$$
\limsup_{n\to \infty} \frac{1}{n} \log |\tr A^n(x)| = \lambda^{+}(A).
$$
By construction, the matrix $H(x)=B(Tx)^{-1}A(x)B(x)$ is diagonal.
We have $A^n(x) = B(T^n x) H^n(x) B(x)^{-1}$ and, in particular,
$\lim\frac{1}{n}\log\|H^n(x)\|=\lambda^{+}(A)$.
Write $B(x)^{-1} B(T^n x) = \left( b_{ij}(n,x) \right)_{i,j=1,2}$.
For a.e. $x$, we know that there are infinitely many $n\in \NN$ such that
$$
\left| b_{11} (n,x) -1 \right|, \,
\left| b_{22} (n,x) -1 \right|< \frac{1}{2}.
$$
The matrices $A^n(x)$ and $B(x)^{-1} B(T^n x) H^n(x)$ have the same trace,
so
\begin{align*}
|\tr A^n(x)| &=
\left| \tr \left( B(x)^{-1} B(T^n x) H^n(x) \right) \right|       \\
&= \big| \left[ b_{11}(n,x) \cdot \|H^n(x)\| +
 b_{22}(n,x) \cdot \|H^n(x)\|^{-1} \right] \big|                  \\
&> \frac{1}{2} \|H^n(x)\| - \frac{3}{2} \|H^n(x)\|^{-1}.        \\
\end{align*}
The result follows.
\qed


\end{document}